\documentclass[11pt,reqno]{amsart}
\usepackage{enumerate, latexsym, amsmath, amsfonts, amssymb, amsthm, color}

\def\Z{\Bbb Z}

\def\bg{\bigg}
\def\({\bg(}
\def\){\bg)}

\def\f{\frac}

\theoremstyle{plain}

\theoremstyle{definition}

\theoremstyle{remark}

\allowdisplaybreaks

 \vspace{4mm}

\begin{document}

\hbox{Amer. Math. Monthly 127(2020), no.\,9, 847--849.}
\medskip

\title
[{Each positive rational number has the form $\varphi(m^2)/{\varphi(n^2)}$}]
{Each positive rational number \\has the form $\varphi(m^2)/{\varphi(n^2)}$}

\author[D. Krachun and Z.-W. Sun] {Dmitry Krachun and Zhi-Wei Sun}

\address {(Dmitry Krachun) St. Petersburg Department of Steklov Mathematical Institute of Russian Academy of Sciences, Fontanka 27, 191023, St. Petersburg, Russia}
\email{dmitrykrachun@gmail.com}

\address{(Zhi-Wei Sun, corresponding author) Department of Mathematics, Nanjing
University, Nanjing 210093, People's Republic of China}
\email{zwsun@nju.edu.cn}

\keywords{Euler's totient function, representation of rational numbers.
%\newline \indent The work is supported by the NSFC-RFBR Cooperation and Exchange Program (grants NSFC 11811530072 and RFBR 18-51-53020-GFEN-a). The second author is also supported
%by the Natural Science Foundation of China (grant 11571162).
}

\begin{abstract}
In this note we show that each positive rational number can be written as $\varphi(m^2)/\varphi(n^2)$,
where $\varphi$ is Euler's totient function and $m$ and $n$ are positive integers.
\end{abstract}
\maketitle

%\section{Introduction}
\setcounter{lemma}{0}
\setcounter{theorem}{0}
\setcounter{corollary}{0}
\setcounter{remark}{0}
\setcounter{equation}{0}

Let $\varphi$ be Euler's totient function. For distinct primes $p_1,\ldots,p_k$
and positive integers $a_1,\ldots,a_k$, it is well known that
$$\varphi(p_1^{a_1}\cdots p_k^{a_k})=\prod_{i=1}^kp_i^{a_i-1}(p_i-1).$$
(See, e.g., \cite[p.\,20]{1}.) Thus, if $n$ has the prime factorization $\prod_{i=1}^kp_i^{a_i}$
(where $p_1,\ldots,p_k$ are distinct primes and $a_1,\ldots,a_k$ are positive integers), then
$$\varphi(n^2)=\prod_{i=1}^kp_i^{2a_i-1}(p_i-1)=n\varphi(n).$$
For positive integers $m$ and $n$ with $\varphi(m^2)=\varphi(n^2)$, we have $m=n$ by comparing
the prime factorizations of $\varphi(m^2)$ and $\varphi(n^2)$. The sequence $\varphi(n^2)\ (n=1,2,3,\ldots)$ is available from \cite{2}.

Section 4 of \cite{3} contains many challenging conjectures on representations of positive rational numbers. For example, Sun \cite[Conjecture 4.4]{3} conjectured that any positive rational number can be written as $m/n$, where $m$ and $n$ are positive integers such that the sum of the $m$th prime and the $n$th prime is a square. Motivated by this, in this note we establish the following new result.

\medskip
\noindent{\bf Theorem 1}. {\it Any positive rational number can be written as $\varphi(m^2)/\varphi(n^2)$,
where $m$ and $n$ are positive integers.}
\medskip

\begin{proof} We claim a stronger result: If $p_1<p_2<\dots < p_k$ are distinct primes and $a_1,\ldots,a_k$ are integers, then there are positive integers $m$ and $n$ with $mn$ not divisible by any prime greater than $p_k$ such that $p_1^{a_1}\cdots p_k^{a_k}=\varphi(m^2)/\varphi(n^2)$.

We prove the claim by induction on $p_k$.

The base of the induction is $p_k=2$. For any $a\in\Z$, clearly
$$2^{2a}=\f{2^{2(a+b)-1}}{2^{2b-1}}=\frac{\varphi(2^{2(a+b)})}{\varphi(2^{2b})}$$
 for each integer $b>|a|$, and
 $$2^{2a+1}=\begin{cases}\varphi(2^{2(a+1)})/\varphi(1^2)&\text{if}\ a\ge0,
 \\\varphi(1^2)/\varphi(2^{-2a})&\text{if}\ a<0.\end{cases}$$

Now let $q$ be an odd prime and assume that the claim holds whenever $p_k<q$.
Let $q_1<\cdots<q_k=q$ be distinct primes
and let $r=\prod_{i=1}^kq_i^{a_i}$ with $a_1,\ldots,a_k\in\Z$. Set
$r_0=r/q_k^{a_k}$ if $2\mid a_k$, and $r_0=r/((q_k-1)q_k^{a_k})$ if $2\nmid a_k$.
Clearly, all the primes in the factorization of $r_0$ are smaller than $q_k=q$.
By the induction hypothesis, there are positive integers $m_0$ and $n_0$ with $m_0n_0$
not divisible by any prime $p\ge q_{k}$ such that
$$\frac{\varphi(m_0^2)}{\varphi(n_0^2)}=r_0.$$
Obviously, we may take $m_0=n_0=1$ if $r_0=1$.

{\it Case} 1. $2\mid a_k$.

In this case, we take positive integers $b$ and $c$ with $b-c=a_k/2$, and set $m=m_0q^b$ and $n=n_0q^c$.
Then
$$\frac{\varphi(m^2)}{\varphi(n^2)}=\frac{\varphi(m_0^2)}{\varphi(n_0^2)}\times \frac{q^{2b-1}(q-1)}{q^{2c-1}(q-1)}=r_0 q^{2(b-c)}=\prod_{i=1}^k q_i^{a_i}=r.
$$

{\it Case} 2. $2\nmid a_k$.

When $a_k>0$, for $m=m_0q^{(a_k+1)/2}$ and $n=n_0$, we have
$$\frac{\varphi(m^2)}{\varphi(n^2)}=\frac{\varphi(m_0^2)}{\varphi(n_0^2)}\times {q^{a_k}(q-1)}=r_0 q^{a_k}(q-1)=\prod_{i=1}^k q_i^{a_i}=r.
$$
If $a_k<0$, then there are positive integers $m$ and $n$ with $mn$ not divisible by any prime greater than $q_k$ such that $\prod_{i=1}^kq_i^{-a_i}=\varphi(n^2)/\varphi(m^2)$ and hence
$\prod_{i=1}^kq_i^{a_i}=\varphi(m^2)/\varphi(n^2)$.
\medskip

In view of the above, the claim holds and hence so does the theorem.
\end{proof}

\noindent {\bf Examples}. We have
$$\frac{19}{47}=\frac{19\times 19673280}{47\times 19673280}=\frac{\varphi(39330^2)}{\varphi(55836^2)}$$
with
$$39330=2\times3^2\times5\times19\times 23
\ \ \text{and}\ \ 55836=2^2\times3^3\times11\times 47.$$
Also,
$$\frac{47}{58}=\frac{47\times 1700160}{58\times 1700160}=\frac{\varphi(14476^2)}{\varphi(20010^2)}$$
with
$$14476=2^2\times7\times11\times 47
\ \ \text{and}\ \ 20010=2\times3\times5\times23\times 29.$$
\medskip

\noindent {\bf Acknowledgment.} The authors would like to thank the two referees for their helpful comments. The work is supported by the NSFC-RFBR Cooperation and Exchange Program (grants NSFC 11811530072 and RFBR 18-51-53020-GFEN-a). The second author is also supported by the National Natural Science Foundation of China (grant no. 11971222).

\setcounter{conjecture}{0} \end{document}